\documentclass[10pt, twoside]{siamltex}

\usepackage{amsmath,amssymb,epic,longtable,curves}

\newfont{\bb}{msbm10}

\begin{document}
\bibliographystyle{plain}
\newtheorem{remark}[theorem]{Remark}
\newtheorem{example}[theorem]{Example}
\newtheorem{prop}[theorem]{Proposition}
\newtheorem{conjecture}[theorem]{Conjecture}

\def\DJ{{\hbox{D\kern-.8em\raise.15ex\hbox{--}\kern.35em}}}
\def\DJo{\DJ okovi\'c}
\def\NSERC{The second author by the NSERC Grant A-5285.}

\def\vf{{\varphi}}
\def\ve{{\varepsilon}}
\def\bR{{\mathbf R}}
\def\bC{{\mathbf C}}
\def\bH{{\mathbf H}}
\def\bZ{{\mathbf Z}}
\def\ad{{\rm ad\;}}
\def\U{{\mbox{\rm U}}}
\def\SO{{\mbox{\rm SO}}}
\def\O{{\mbox{\rm O}}}
\def\Sp{{\mbox{\rm Sp}}}

\title
{Proof of Atiyah's conjecture for two special
types of configurations}

\author{Dragomir \v Z. \DJ okovi\'c}\thanks
{Department of Pure Mathematics, University of Waterloo,
Waterloo, Ontario, N2L 3G1, Canada.
E-mail: djokovic@uwaterloo.ca.
\newline\indent
The author was supported in part by
the NSERC Grant A-5285.}

\markboth{D.\v Z.\ Djokovi\'c}
{Atiyah's conjecture}

\maketitle

\begin{abstract}
To an ordered $N$-tuple $(x_1,\ldots,x_N)$ of distinct points in
$\bR^3$, Atiyah \cite{MA1,MA2} has associated an ordered
$N$-tuple of homogeneous polynomials
$(p_1,\ldots,p_N)$ in $\bC[x,y]$ of degree $N-1$,
each $p_i$ determined only up to a scalar factor.
He has conjectured that these polynomials are linearly
independent. We show that his conjecture is true for two
special configurations of $N$ points.
Moreover we show that, for one of these configurations,
the stronger conjecture \cite[Conjecture 2]{AS}
is valid.
\end{abstract}

\begin{keywords}{Atiyah's conjecture, the Hopf map,
configuration of $N$ points in $\bR^3$, projective line $\bf{PC}^1$.}
\end{keywords}

\begin{AMS} Primary 51M04, 51M16, Secondary 70G25 \end{AMS}

\date{}

\section{Two conjectures}

Let $(x_1,\ldots,x_N)$ be an ordered $N$-tuple of distinct
points in $\bR^3$. Each ordered pair
$(x_i,x_j)$ with $i\ne j$ determines a point
$$
\frac{x_j-x_i}{|x_j-x_i|}
$$
on the unit sphere $S^2\subset\bR^3$. Identify $S^2$ with the
complex projective line ${\bf PC}^1$ by using a stereographic
projection. Hence one obtains a point 
$(u_{ij},v_{ij})\in{\bf PC}^1$ and a nonzero linear form
$l_{ij}=u_{ij}x+v_{ij}y\in\bC[x,y]$.
Define homogeneous polynomials $p_i\in\bC[x,y]$ of
degree $N-1$ by
\begin{equation} \label{poli}
p_i=\prod_{j\ne i}l_{ij}(x,y),\quad i=1,\ldots,N.
\end{equation}

\begin{conjecture} \label{conj-1} (Atiyah \cite{MA2})
The polynomials $p_1,\ldots,p_N$ are linearly independent.
\end{conjecture}

Atiyah \cite{MA1,MA2} has observed that his conjecture is true
if the points $x_1,\ldots,x_N$ are collinear. He has also
verified the conjecture for $N=3$. The case $N=4$ has
been verified by Eastwood and Norbury \cite{EN}. For
additional information on the conjecture (further
conjectures, generalizations, and numerical evidence)
see \cite{MA2,AS}.

In order to state the second conjecture, one has to be more explicit.
Identify $\bR^3$ with $\bR\times\bC$ and denote the origin by $O$.
Following Eastwood and Norbury \cite{EN}, we make use of the Hopf map
$h:\bC^2\setminus\{O\}\to(\bR\times\bC)\setminus\{O\}$ defined by:
$$
h(z,w)=((|z|^2-|w|^2)/2,z\bar{w}).
$$
This map is surjective and its fibers are the circles
$\{(zu,wu):u\in S^1\}$, where $S^1$ is the unit circle in $\bC$.
If $h(z,w)=(a,v)$, we say that $(z,w)$ is a {\em lift} of $(a,v)$.
For instance, we can take
$$
\lambda^{-1/2}(\lambda,\bar{v}),\quad
\lambda=a+\sqrt{a^2+|v|^2},
$$
as the lift of $(a,v)$.

Assume that our points are $x_i=(a_i,z_i)$.
For the sake of simplicity assume that if $i<j$ and
$z_i=z_j$ then $a_i<a_j$.
As the lift of the vector $x_j-x_i$, $i<j$, we choose
$$
\frac{1}{\sqrt{\lambda_{ij}}}\left(\lambda_{ij},
\bar{z}_j-\bar{z}_i\right),
$$
where
$$
\lambda_{ij}=a_j-a_i+\sqrt{(a_j-a_i)^2+|z_j-z_i|^2}.
$$
According to the recipe in \cite{MA2,AS,EN}, we always
use the lift $(-\bar{w},\bar{z})$ for the vector $x_i-x_j$
if $(z,w)$ has been chosen as the lift of $x_j-x_i$.
Hence we introduce the linear forms
\begin{eqnarray*}
l_{ij}(x,y) &=& \lambda_{ij}x+(\bar{z}_j-\bar{z}_i)y,\quad i<j; \\
l_{ij}(x,y) &=& (z_j-z_i)x+\lambda_{ji}y,\quad i>j.
\end{eqnarray*}

Define $P$ to be the $N\times N$ coefficient matrix of the binary forms
$p_i(x,y)$ defined by (\ref{poli}) using the above $l_{ij}$'s.
The second conjecture that we are interested in can now be
formulated as follows.

\begin{conjecture} \label{conj-2}
(Atiyah and Sutcliffe \cite[Conjecture 2]{AS}, see also \cite{EN}).
If $r_{ij}=|x_j-x_i|$, then
$$
|\det(P)|\ge\prod_{i<j}(2\lambda_{ij}r_{ij}).
$$
\end{conjecture}

As $2\lambda_{ij}r_{ij}=\lambda_{ij}^2+|z_j-z_i|^2$, this conjecture can
be rewritten as
\begin{equation} \label{conj2}
|\det(P)|\ge\prod_{i<j}\left(\lambda_{ij}^2+|z_j-z_i|^2\right).
\end{equation}

Obviously, this conjecture is stronger than Conjecture \ref{conj-1}.

\section{Two special cases of Atiyah's conjecture}

We shall prove Atiyah's conjecture in the following two cases:

\begin{itemize}
\item[(A)] $N-1$ of the points $x_1,\ldots,x_N$ are collinear.
\item[(B)] $N-2$ of the points $x_1,\ldots,x_N$ are on a line
$L$ and the line segment joining the remaining two points has
its midpoint on $L$ and is perpendicular to $L$.
\end{itemize}

Let $L$ and $M$ be two perpendicular lines in $\bR^3$
intersecting at the origin, $O$.
Let $N=m+n$ and assume that the points
$x_1,\ldots,x_m$ are on $L$ and $x_{m+1},\ldots,x_N$ are
on $M$ but not on $L$.
Set $y_j=x_{m+j}$ for $j=1,\ldots,n$.

Without any loss
of generality, we may assume that $L=\bR\times\{0\}$
and $M=\{0\}\times\bR$.
Write $x_i=(a_i,0)$ for $i=1,\ldots,m$ and $y_j=(0,b_j)$
for $j=1,\ldots,n$.
We may also assume that $a_1<a_2<\cdots<a_m$
and $b_1<b_2<\cdots<b_n$.

The lifts of the nonzero vectors $x_j-x_i$,
$i,j\in\{1,\ldots,N\}$ are given in Table 1,
where we have set
$$
\lambda_{ij}=a_i+\sqrt{a_i^2+b_j^2}.
$$

\vspace{3mm}

\begin{center}

{\bf Table 1: The lifts of the vectors $x_j-x_i$ }
$$
\begin{array}{|l|c|r|c|}
\hline
{\rm Vectors} & {\rm Index~ restrictions} &
{\rm Lifts}\quad\quad\quad & {\rm Linear~ forms} \\
\hline
x_r-x_i & 1\le i<r\le m & \left(2(a_r-a_i)\right)^{1/2}(1,0) &
2(a_r-a_i)x \\
x_i-x_r & 1\le i<r\le m & \left(2(a_r-a_i)\right)^{1/2}(0,1) &
2(a_r-a_i)y \\
y_s-y_j & 1\le j<s\le n & (b_s-b_j)^{1/2}(1,1) & (b_s-b_j)(y+x) \\
y_j-y_s & 1\le j<s\le n & (b_s-b_j)^{1/2}(-1,1) & (b_s-b_j)(y-x) \\
x_i-y_j & 1\le i\le m,~ 1\le j\le n &
\lambda_{ij}^{-1/2}(\lambda_{ij},-b_j) & \lambda_{ij}x-b_jy \\
y_j-x_i & 1\le i\le m,~ 1\le j\le n &
\lambda_{ij}^{-1/2}(b_j,\lambda_{ij}) & b_jx+\lambda_{ij}y \\
\hline
\end{array}
$$
\end{center}
\vspace{3mm}

The associated polynomials $p_i$ (up to scalar factors)
are given by:
\begin{eqnarray}
&& p_i(x,y)=x^{m-i}y^{i-1}\prod_{j=1}^n
(b_jx+\lambda_{ij}y),\quad 1\le i\le m; \label{jed-m} \\
&& p_{m+j}(x,y)=(y+x)^{n-j}(y-x)^{j-1}\prod_{i=1}^m
(\lambda_{ij}x-b_jy), \quad 1\le j\le n. \label{jed-n}
\end{eqnarray}

\begin{theorem} \label{T-A}
Conjecture \ref{conj-1} is valid under the hypothesis (A).
\end{theorem}

\begin{proof}
In this case we have $n=1$.
Without any loss of generality we may assume that $b_1=-1$.
After dehomogenizing the polynomials
$p_i$ (or $-p_i$) by setting $x=1$, we obtain the polynomials:
\begin{eqnarray*}
&& y^{i-1}(1-\lambda_iy),\quad 1\le i\le m; \\
&& \prod_{i=1}^m (y+\lambda_i),
\end{eqnarray*}
where $\lambda_i=\lambda_{i1}>0$.
The coefficient matrix of these polynomials is:
$$
\begin{pmatrix}
1 & -\lambda_1 & 0 & 0 & \ldots & 0 & 0 \cr
0 & 1 & -\lambda_2 & 0 &        & 0 & 0 \cr
0 & 0 & 1 & -\lambda_3 &        & 0 & 0 \cr
\vdots & & & & & & \cr
0 & 0 & 0 & 0		& 	& 1 & -\lambda_m \cr
E_m & E_{m-1} & E_{m-2} & E_{m-3} & & E_1 & 1 \cr
\end{pmatrix}
$$
where $E_k$ is the $k$-th elementary symmetric function of
$\lambda_1,\ldots,\lambda_m$. Its determinant,
$$
1+\lambda_mE_1+\lambda_{m-1}\lambda_mE_2+\cdots+
\lambda_1\lambda_2\cdots\lambda_mE_m,
$$
is positive.
\end{proof}

\begin{theorem} \label{T-B}
Conjecture \ref{conj-1} is valid under the hypothesis (B).
\end{theorem}

\begin{proof}
In this case $n=2$ and $b_1+b_2=0$.
Without any loss of generality we may assume that $b_1=-1$.
After dehomogenizing the polynomials
$p_i$ (or $-p_i$) by setting $x=1$, we obtain the polynomials:
\begin{eqnarray*}
&& y^{i-1}(1-\lambda_i^2y^2),\quad 1\le i\le m; \\
&& (y+1)\prod_{i=1}^m (y+\lambda_i), \\
&& (y-1)\prod_{i=1}^m (y-\lambda_i),
\end{eqnarray*}
where $\lambda_i=\lambda_{i1}>0$.
The coefficient matrix of these polynomials is:
$$
\begin{pmatrix}
1 & 0 & -\lambda_1^2 & 0 & \ldots & 0 & 0 & 0 \cr
0 & 1 & 0 & -\lambda_2^2 &        & 0 & 0 & 0 \cr
\vdots & & & & & & & \cr
0 & 0 & 0 & & & 1 & 0 & -\lambda_m^2 \cr
\tilde{E}_{m+1} & \tilde{E}_m & \tilde{E}_{m-1} & & & \tilde{E}_2
& \tilde{E}_1 & 1 \cr
(-1)^{m+1}\tilde{E}_{m+1} & (-1)^m\tilde{E}_m & (-1)^{m-1}\tilde{E}_{m-1} &
 & & \tilde{E}_2 & -\tilde{E}_1 & 1 \cr
\end{pmatrix}
$$
where $\tilde{E}_k$ is the $k$-th elementary symmetric function of
$1,\lambda_1,\ldots,\lambda_m$. Its determinant is $2pq$ where
\begin{eqnarray*}
p &=& 1+\lambda_m^2\tilde{E}_2+\lambda_{m-2}^2\lambda_m^2\tilde{E}_4+\cdots,\\
q &=& \tilde{E}_1+\lambda_{m-1}^2\tilde{E}_3+
\lambda_{m-3}^2\lambda_{m-1}^2 \tilde{E}_5+\cdots, 
\end{eqnarray*}
and so it is positive.
\end{proof}

\section{Atiyah and Sutcliffe conjecture is valid in case (A)}

In the general setup of the previous
section, the Conjecture \ref{conj-2} asserts that
\begin{equation} \label{opsta}
|\det(P)|\ge 2^{\binom n2}\prod_{i,j}\left(
\lambda_{ij}^2+b_j^2 \right).
\end{equation}
where $P$ is the coefficient matrix (of order $N=m+n$) of 
the polynomials (\ref{jed-m}) and (\ref{jed-n}).

In case (A) this inequality takes the form
\begin{equation} \label{nej-A}
1+\lambda_mE_1+\lambda_{m-1}\lambda_mE_2+\cdots+
\lambda_1\lambda_2\cdots\lambda_mE_m\ge
\prod_{i=1}^m (1+\lambda_i^2),
\end{equation}
where, as in the proof of Theorem \ref{T-A}, we assume that
$b_1=-1$ and $E_k$ denotes the $k$-th elementary symmetric
function of $\lambda_1,\ldots,\lambda_m$. Thus we have
$$
\lambda_i=a_i+\sqrt{1+a_i^2}>0
$$
and
\begin{equation} \label{lambda}
\lambda_1<\lambda_2<\cdots<\lambda_m.
\end{equation}
Let $E^{(2)}_k$ denote the $k$-th elementary symmetric
function of $\lambda_1^2,\ldots,\lambda_m^2$.
In view of (\ref{lambda}), we have
$$
\lambda_{m-k+1}\lambda_{m-k+2}\cdots\lambda_mE_k\ge E^{(2)}_k,
\quad 0\le k\le m.
$$
The inequality (\ref{nej-A}) is a consequence of the inequalities
just written since
$$
\prod_{i=1}^m(1+\lambda_i^2)=\sum_{k=0}^m E^{(2)}_k.
$$
Hence we have the following result.

\begin{theorem}
Conjecture \ref{conj-2} is valid in case (A).
\end{theorem}

In case (B) the inequality (\ref{opsta}) takes the form:
\begin{eqnarray*}
&& \left(1+\lambda_m^2\tilde{E}_2+\lambda_{m-2}^2\lambda_m^2\tilde{E}_4
+\cdots\right)\cdot\left(\tilde{E}_1+\lambda_{m-1}^2\tilde{E}_3+
\lambda_{m-3}^2\lambda_{m-1}^2 \tilde{E}_5+\cdots \right) \\
&& \quad\ge\prod_{i=1}^m\left(1+\lambda_i^2\right)^2,
\end{eqnarray*}
where $\tilde{E}_k$ are as in the proof of Theorem \ref{T-B}.

It is easy to verify that this inequality holds for $m=1$,
but we were not able to prove it in general.
If we set all $\lambda_i=\lambda>0$,
then the above inequality specializes to
\begin{eqnarray*}
&& \left[(1+\lambda^2)^m+\sum_{k\ge0}\binom{m}{2k+1}
(\lambda^{4k+3}-\lambda^{4k+2})\right]\cdot \\
&& \left[(1+\lambda^2)^m-\sum_{k\ge0}\binom{m}{2k+1}
(\lambda^{4k+2}-\lambda^{4k+1})\right]\ge(1+\lambda^2)^{2m}.
\end{eqnarray*}
Since
$$
\sum_{k\ge0}\binom{m}{2k+1}(\lambda^{4k+3}-\lambda^{4k+2})
=\frac{1}{2}(\lambda-1)\left[(1+\lambda^2)^m-(1-\lambda^2)^m\right],
$$
it is easy to verify that the specialized inequality is valid.



\begin{thebibliography}{99}

\bibitem{MA1}
M. Atiyah.
\newblock The geometry of classical particles.
\newblock {\em Surveys in Differential Geometry}
(International Press) {\bf 7} (2001).

\bibitem{MA2}
M. Atiyah.
\newblock Configurations of points.
\newblock {\em Phil. Trans. R. Soc. Lond. A} 
{\bf 359} (2001), 1375--1387.

\bibitem{AS}
M. Atiyah and P. Sutcliffe.
\newblock The geometry of point particles.
\newblock {\bf arXiv:hep-th/0105179 v2}, 15 Oct 2001.

\bibitem{EN}
M. Eastwood and P. Norbury.
\newblock A proof of Atiyah's conjecture on configurations
of four points in Euclidean three-space.
\newblock {\em Geometry \& Topology} {\bf 5} (2001), 885--893.



\end{thebibliography}
\end{document}